\documentclass[10pt, titlepage]{amsart}

\usepackage{ae} % or {zefonts}
\usepackage[T1]{fontenc}
\usepackage[cp1250]{inputenc}
\usepackage{amsmath}
\usepackage{amssymb, amsfonts,amscd,verbatim}
\usepackage[dvips]{graphicx}
\usepackage{latexsym}
\usepackage{indentfirst}
\usepackage{latexsym}

\usepackage{graphicx}
\usepackage{verbatim}   % useful for program listings
\usepackage{color}      % use if color is used in text
\usepackage{subfigure}  % use for side-by-side figures

\usepackage{amsmath}    % need for subequations

\theoremstyle{plain}
\newtheorem{Prop}{Proposition}[section]
\newtheorem{Thm}[Prop]{Theorem}
\newtheorem{Cor}[Prop]{Corollary}

\newtheorem{Lem}[Prop]{Lemma}

\theoremstyle{definition}

\theoremstyle{remark}

\newtheorem{Question}[Prop]{\bf Question}

\def\int{\mathop{\roman{int}}}

\def\1{^{-1}}

\def\Z{{\mathbf Z}}

\def\R{{\mathbf R}}
\def\NN{{\mathbf N}}
\def\RR{{\mathbf R}}
\def\HH{{\mathbf H}}

\def\dim{\text{dim}}

\def\UUU{{\mathcal U}}

\def\asdim{\mathrm{asdim}}

\def\AC{\mathrm{Cone}}
\def\dim{\mathrm{dim}}

\def\dokaz{{\bf Proof. }}
\def\edokaz{\hfill $\blacksquare$}

\numberwithin{equation}{section}

\input pstricks.tex
\input xy
\xyoption{all}

\begin{document}
%\fontsize{20}{34pt}\selectfont
\title[
Asymptotic cones and Assouad-Nagata dimension
]%
   {Asymptotic cones and Assouad-Nagata dimension}

\author{J.~Dydak}
\address{University of Tennessee, Knoxville, TN 37996, USA}
\email{dydak@math.utk.edu}

\author{J.~Higes}
\address{Departamento de Geometr\'{\i}a y Topolog\'{\i}a,
Facultad de CC.Matem\'aticas. Universidad Complutense de Madrid.
Madrid, 28040 Spain} \email{josemhiges@yahoo.es}
\date{ \today
}


\thanks{The first-named author was partially supported
by Grant No.  2004047  from the United States-Israel Binational Science
Foundation (BSF),  Jerusalem, Israel. The second named author is supported by Grant AP2004-2494 from the Ministerio de Educaci\' on y Ciencia, Spain.}

\begin{abstract}

We prove the dimension
of any asymptotic cone over a metric space $(X,\rho)$
does not exceed the asymptotic Assouad-Nagata dimension 
$\asdim_{AN}(X)$ of $X$.
This improves a result of Dranishnikov and Smith \cite{DRS}
who showed $\dim(Y)\leq \asdim_{AN}(X)$ for all
separable subsets $Y$ of special asymptotic cones
$\AC_\omega(X)$, where $\omega$ is an exponential
ultrafilter on natural numbers.
\par We also show that Assouad-Nagata dimension of the discrete
Heisenberg group equals its asymptotic dimension.
\end{abstract}
\maketitle

\medskip
%Printed on \today.
\medskip
\tableofcontents

\section{Introduction}

Given a metric space $(X,\rho_X)$ and given a sequence $d=\{d_n\}_{n=1}^\infty$
of real numbers diverging to infinity, one constructs
{\it asymptotic cone} $\AC_\omega(X,c,d)$ for every non-principal
ultrafilter $\omega$ in natural numbers $\NN$ as follows:
\par First one considers the pseudo-metric space of
all functions $x\colon N\to X$ such that $\{\frac{\rho_X(x(n),c(n))}{d_n}\}$
is bounded for a given function $c\colon N\to X$ (that function serves as the base point
of the asymptotic cone).
The distance from $x$ to $y$ is the $\omega$-limit 
$\lim_\omega \frac{\rho_X(x(n),y(n))}{d_n}$ of
$\{\frac{\rho_X(x(n),y(n))}{d_n}\}$ (see \cite{Behr} or \cite{DruSap} for details). 
By identifying $x$ and $y$ whose distance is $0$ one gets a metric space
$\AC_\omega(X,c,d)$.
\par One can generalize the definition of asymptotic cones
for non-principal ultrafilters $\omega$ over any infinite set $I$ (see \cite{Kap}).
In that case one requires that $d=\{d_i\}_{i\in I}$ has the property
that $\omega$-limit of $\{\frac{1}{d_i}\}_{i\in I}$ is $0$.
As all our proofs are valid for that case, we will work in that generality.
\par An important problem in geometric group theory is relating
properties of finitely generated groups $G$ to topological
properties of their asymptotic cones (i.e., asymptotic cones of their Cayley graphs):
\begin{itemize}
\item A finitely generated group is virtually Abelian if and only if its asymptotic cones are isometric to the space $\RR^n$  \cite{GroNilp} \cite{Pansu}.

\item A finitely generated group is virtually nilpotent if and only if its asymptotic cones are locally compact
\cite{GroNilp} \cite{Dries} \cite{Drutu}.

\item A finitely generated group is hyperbolic if and only if its asymptotic cones are $\RR$-trees \cite{GroAsymInvar}.

\item A finitely generated group $G$ is relatively hyperbolic with respect to finitely
generated subgroups $H_1,\ldots, H_n$ if and only if every asymptotic cone of $G$ is tree-graded with
respect to $\omega$-limits of sequences of cosets of the subgroups $H_i$ \cite{DruSap}.
\end{itemize}

It would be of interest to detect the property of a finitely generated group $G$
equivalent to finite-dimensionality of its asymptotic cones.

The first time dimension of asymptotic cones was discussed
occurred in \cite{GroNilp}, where M.Gromov proved finite-dimensionality
of asymptotic cones of a finitely generated group of polynomial growth
by demonstrating that their Hausdorff dimension is finite.

J.Burillo \cite{Bur} showed $1$-dimensionality of asymptotic cones of Baumslag-Solitar groups
 and used it to demonstrate that the fundamental group
of those cones is not free and uncountable.

More recently, J.Behrstock and Y.Minsky \cite{BehrMinsky}
proved the following:

\begin{Thm}[Dimension Theorem of Behrstock-Minsky]\label{BMThm}
 The maximal topological dimension of a locally compact
subset of the asymptotic cone of a mapping class group is equal to
the maximal rank of an abelian subgroup.
\end{Thm}

The aim of this note is to relate asymptotic Assouad-Nagata dimension of a finitely
generated group $G$ to the dimension of its asymptotic cones. More generally,
we prove
$\dim(\AC_\omega(X,c,d))\leq \asdim_{AN}(X)$ for all 
metric spaces $(X,\rho)$ and all non-principal ultrafilters $\omega$.
Its direct predecessor is the following result from \cite{DRS}:

\begin{Thm}[Dranishnikov-Smith]\label{DRSCor}
$\dim(Y)\leq \asdim_{AN}(X)$ for every separable subset
$Y\subset \AC_\omega(X)$ if $X$ is a proper metric space and $\omega$
is an exponential ultrafilter in $\NN$.
\end{Thm}

Notice that asymptotic cones in \ref{DRSCor} are based on
the sequence $d_n=n$.

The way \ref{DRSCor} was proved in \cite{DRS} is
by constructing an injective map $\AC_\omega(X)\setminus [x_0]\to \nu_L(X)$,
where $\nu_L(X)$ is the {sublinear Higson corona} of $X$
(see Theorem 5.1 of \cite{DRS}). 
Since $\asdim_{AN}(X)\ge \dim(\nu_L(X))$ for every proper metric space $(X,d)$
(see Theorem 3.10 in \cite{DRS}),  \ref{DRSCor} follows.

It is not clear if the methods of this paper allow to improve
\ref{BMThm} and show that the dimension of
 asymptotic cones of mapping class groups equals
 the maximal rank of an abelian subgroup.
 We learned from Jason Behrstock that in order to  remove the \lq locally compact\rq\  hypothesis 
in \ref{BMThm} all one needs to show is that the dimension
of asymptotic cones of mapping class groups is finite.
That leads to the question of Assouad-Nagata dimension of mapping class
groups being finite. That question remains open. Even
the finiteness of asymptotic dimension of mapping class groups
remains open (except in genus less than 3).

\par We are grateful to Jason Behrstock and Kolya Brodskiy for valuable comments.
 
\section{Assouad-Nagata dimension}

A metric space $(X,\rho)$ is of {\it Assouad-Nagata dimension}
at most $n$ (notation: $\dim_{AN}(X)\leq n$)
if there is a constant $c > 0$ such that for any $r > 0$
there is a cover $\UUU_r$ of $X$ whose elements have
diameter at most $c\cdot r$ and every $r$-ball $B(x,r)$
intersects at most $n+1$ elements of $\UUU_r$. Equivalently (see \cite{BDLM}),
there is a constant $K > 0$ such that for any $r > 0$
$X$ can be expressed as $X=X_0\cup\ldots\cup X_n$
so that all $r$-components of $X_i$ are $K\cdot r$-bounded (two points
are in the same $r$-component of a metric space if they can be connected
by a finite chain of points so that each link is of size less than $r$).

\par If the existence of $\UUU_r$ is guaranteed only for sufficiently
large $r$, then we say {\it the asymptotic Assouad-Nagata dimension}
of $X$ is at most $n$ (notation: $\asdim_{AN}(X)\leq n$).

\par Of basic importance to us is the following result:

\begin{Thm}[Lang-Schlichenmaier \cite{LangSch}]\label{LSThm}
If $f\colon (X,\rho_X)\to (Y,\rho_Y)$
is a quasi-symmetric embedding of metric spaces,
then $\dim_{AN}(X)\leq \dim_{AN}(Y)$.
\end{Thm}

We will apply \ref{LSThm} in the case of bi-Lipschitz homeomorphisms
and snowflaking. Recall that a {\it snowflaked version} of a metric space
$(X,\rho)$ is the space $(X,\rho^p)$ for some $0 < p < 1$
(see \cite{Hei}, p.98).

$(X,\rho)$ has {\it $n$-dimensional Nagata property}
if for every configuration $y_1,\ldots,y_{n+2}$
of $n+2$ points in $X$ the conditions
$\rho(y_i,B(x,r/2))< r$ for some $x\in X$ and all $i\leq n+2$
imply existence of $i\ne j$ such that $\rho(y_i,y_j) < r$.

\begin{Thm}\label{ANDimAndPnProperty}
Let $(X,\rho_X)$ be a metric space and $n\ge 0$.
$\dim_{AN}(X)\leq n$ if and only if 
 $(X,\rho_X)$ is bi-Lipschitz
homeomorphic to a metric space $(Y,\rho_Y)$ such that some snowflaked version of $(Y,\rho_Y)$
has $n$-dimensional Nagata property.
\end{Thm}
\dokaz In one direction it follows from \ref{LSThm} and the fact
that $\dim_{AN}(Y)\leq n$ if $Y$ has $n$-dimensional Nagata property.
Indeed, if $r > 0$ and $A$ is a maximal subset of $Y$
with respect to the property $\rho(x,y)\ge r$ for all $x\ne y\in Y$,
then the family of balls $\UUU_r=\{B(y,r)\}_{y\in A}$ covers $Y$
and any ball $B(x,r/2)$ intersects at most $n+1$ elements of $\UUU_r$
(see the proof of Theorem 4.2 in \cite{Dran-Zar} along the same lines).
\par The other direction is proved in Proposition 2.2 of \cite{Assouad}.
\edokaz

\section{Maps with metrically parallel fibers}

J.Burillo \cite{Bur} defined a very interesting class of maps, maps with metrically parallel fibers.
Recall $f\colon (X,\rho_X)\to (Y,\rho_Y)$ is a map
{\it with metrically parallel fibers}
if for every $x\in X$ and $y\in f(X)$ there is $x'\in f^{-1}(y)$
such that $\rho_X(x,x')=\rho_Y(f(x),f(x'))=\rho_Y(f(x),y)$
(see Definition 6 in \cite{Bur}).
It seems maps with metrically parallel fibers provide
a natural ground for applying Hurewicz type theorems from \cite{BDLM}
(see the proof of \ref{GeneralDimForMetricallyParallelFibers}).

\par {\bf Examples of maps with metrically parallel fibers}:
\begin{enumerate}
\item The projection $X\times Y\to Y$ if the metric on $X\times Y$
is either euclidean or the sum of corresponding metrics on $X$ and $Y$.
\item $Im\colon \HH^2\to \RR$, where $\HH^2$ is the hyperbolic plane
in the upper half plane model (see \cite{Bur}). Notice $\RR$
is identified with $(0,\infty)$ equipped with the metric $d(x,y)=|\ln(x)-\ln(y)|$.
\end{enumerate}

We will weaken the condition of existence of $x'$ such that $\rho_X(x,x')=\rho_Y(f(x),f(x'))=\rho_Y(f(x),y)$ to simply $\rho_X(x,f^{-1}(y))\leq \rho_Y(f(x),y)$
 for all $x\in X$ and $y\in f(X)$.
 One reason for the change is that the following natural functions
 easily satisfy the new condition.

\begin{enumerate}
\item The function $x\to d(x,x_0)$ from any $\RR$-tree $T$ to $\RR_+$,
where $x_0$ is a base point of $T$. We assume every point of $T$
lies inside an infinite geodesic. 
\item Any epimorphism $f\colon (G,d_G)\to (H,d_H)$
of finitely generated groups equipped with word metrics induced by
finite sets of generators $S$ in $G$ and $f(S)$ in $H$.
\end{enumerate}

The second reason is the following.
\begin{Prop}\label{ReasonForMetricallyParallelFibers}
If $f\colon (X,\rho_X)\to (Y,\rho_Y)$ is a Lipschitz map, then the following
conditions are equivalent:
\begin{itemize}
\item[a.] The natural function $f^{-1}(y)\to y$
establishes a bi-Lipschitz equivalence between
$f(X)$ and the space of fibers of $f$ equipped with the Hausdorff distance.
\item[b.] There is $\mu > 0$ such that $\rho_X(x,f^{-1}(y))\leq \mu\cdot \rho_Y(f(x),y)$
 for all $x\in X$ and $y\in f(X)$.
\end{itemize}
\end{Prop}
\dokaz a)$\implies$b). Recall the Hausdorff distance $d_H(A,B)$ between two subsets
$A$ and $B$ of $X$ is the supremum of $\rho_X(a,B)$ and $\rho_X(b,A)$
over all $a\in A$ and $b\in B$.
A surjective map $g\colon S\to T$ of metric spaces is a bi-Lipschitz
equivalence if there are constants $\lambda,\mu > 0$
such that $\lambda\cdot d_T(g(a),g(b))\leq d_S(a,b)\leq \mu\cdot d_T(g(a),g(b))$
for all $a,b\in S$.
Thus b) follows easily from a) as $\rho_X(x,f^{-1}(y))\leq d_H(f^{-1}(f(x)),f^{-1}(y))\leq \mu\cdot \rho(f(x),y)$ for all $x\in X$ and $y\in f(X)$.
\par b)$\implies$a). Assume $f$ is $\lambda$-Lipschitz and $y_1,y_2\in f(X)$.
Since $\rho_X(x,f^{-1}(y_2))\leq \mu\cdot \rho(y_1,y_2)$
for all $x\in f^{-1}(y_1)$ and $\rho_X(x,f^{-1}(y_1))\leq \mu\cdot \rho(y_1,y_2)$
for all $x\in f^{-1}(y_2)$, we get $d_H(f^{-1}(y_1), f^{-1}(y_2))\leq \mu\cdot \rho_Y(y_1,y_2)$.
Similarly, as $f$ is $\lambda$-Lipschitz,
one has $\lambda\cdot \rho_X(x,f^{-1}(y_2))\ge \rho_Y(y_1,y_2)$
for all $x\in f^{-1}(y_1)$ resulting in $\frac{\rho_Y(y_1,y_2)}{\lambda}\leq d_H(f^{-1}(y_1), f^{-1}(y_2))$ for all $y_1,y_2\in f(X)$.
\edokaz

\ref{ReasonForMetricallyParallelFibers} provides a global explanation
for the term \lq metrically parallel fibers\rq. \ Also, after rescaling $Y$
one can simply achieve $\mu=1$ in \ref{ReasonForMetricallyParallelFibers},
so that is going to be our assumption.

It was suggested to us by K.Brodskiy that maps as in \ref{ReasonForMetricallyParallelFibers}(b)
ought to be called {\bf metrically open} due to the following:

\begin{Lem}[K.Brodskiy]
Let $f\colon X\to Y$ be a surjective map of metric spaces. The following
properties are equivalent:
\begin{itemize}
\item[1.] There exists $\mu>0$ such that for any $x\in X$ and $y\in Y$,
$$ d_X(x,f^{-1}(y))\le \mu d_Y(f(x),y). $$
\item[2.] There exists $\lambda>0$ such that for any $x\in X$ and any $R>0$,
$$ B(f(x),\lambda R)\subset f(B(x,R)) $$
\end{itemize}
\end{Lem}
\dokaz  In both directions assume relation $\mu=1/\lambda$.
\par 1)$\implies$2). If $d_Y(f(x),f(z)) < \lambda R$ and $f(z)=y$,
then $d_X(x,f^{-1}(y)) < R$, so there is $t\in f^{-1}(y)$ satisfying $d_X(x,t) < R$.
Thus $y=f(t)\in f(B(x,R))$.
\par 2)$\implies$1). Suppose $d_Y(f(x),y) < \lambda R$. That means $y\in B(f(x), \lambda R)$
and there is $z\in B(x,R)$ so that $y=f(z)$.
Thus $d_X(x,f^{-1}(y)) < R$ which proves $ d_X(x,f^{-1}(y))\le \mu d_Y(f(x),y). $
\edokaz

J.Burillo proved in \cite{Bur} (Theorem 16) that $\dim(X)\leq \dim(Y)$
if there is a map $f\colon X\to Y$ with parallel fibers such that
$f^{-1}(y)$ is an ultrametric space for all $y\in Y$. Recall $(A,\rho)$
is {\it ultrametric} if $\rho(x,z)\leq \max(\rho(x,y),\rho(y,z))$
for all triples $x,y,z$ in $A$.
Our analog of that result for Assouad-Nagata dimension has a much simpler
proof (if Lipschitz maps are involved):

\begin{Thm}\label{DimZeroForMetricallyParallelFibers}
Suppose $f\colon (X,\rho_X)\to (Y,\rho_Y)$ is a Lipschitz map with ultrametric fibers.
If $$\rho_X(x,f^{-1}(y))\leq \rho_Y(f(x),y)$$
 for all $x\in X$ and $y\in f(X)$,
then $\dim_{AN}(X)\leq \dim_{AN}(Y)$.
\end{Thm}
\dokaz Let $f$ be $\lambda$-Lipschitz (that means $\rho(f(x),f(y))\leq \lambda\cdot \rho_X(x,y)$
for all $x,y\in X$).
Let us show that if $B\subset Y$ and all $\lambda\cdot r$-components of $B$
are $K$-bounded, then all $r$-components of $f^{-1}(B)$ are
$(4K+r)$-bounded.
Indeed, if $x_1,\ldots , x_n$ is a chain of points such that $f(x_i)\in B$
and $\rho_X(x_i,x_{i+1}) < r$ for all $i< n$, then $f(x_i)$ belong to the same
$\lambda\cdot r$-component of $B$. Hence there is a fiber
$F\subset f^{-1}(B)$ such that $\rho_X(x_i,F) < K$ for all $i\leq n$
and we can pick $y_i\in F$ satisfying $\rho_X(x_i,y_i) < K$ for all $i\leq n$.
Thus all $y_i$ belong to the same $(2K+r)$-component of $F$.
As $F$ is ultrametric, that component is of size less than $2K+r$.
Therefore $\rho_X(x_1,x_n)\leq 2K+\rho_X(y_1,y_n)< 2K+2K+r=4K+r$.
\par
Suppose $\dim_{AN}(Y)=m < \infty$ (the case $\dim_{AN}(Y)=\infty$
is trivial) and let $c > 0$ be the corresponding constant for $Y$.
Given $r > 0$ pick a cover $\{Y_1,\ldots,Y_{m+1}\}$ of $Y$
such that any $\lambda\cdot r$-component of each $Y_i$ has diameter less than $c\cdot \lambda\cdot r$.
Now $r$-components of each set $f^{-1}(Y_i)$ are bounded by
$(4c\lambda+1)\cdot r$ which proves $\dim_{AN}(X)\leq m$.
\edokaz

Notice \ref{DimZeroForMetricallyParallelFibers} implies the well-known fact
$\dim_{AN}(T)\leq 1$ for any $\RR$-tree $T$ using
$f\colon T\to \RR_+$, $f(x)=d(x,x_0)$.

Bell and Dranishnikov \cite{BellDranish} introduced the concept
of
$\asdim\{A_s\}_{s\in S} \le n$ uniformly with respect to $s\in S$.
The analogous concept of {\it $\dim_{AN}(A_s)\leq n$ uniformly
with respect to $s\in S$} means that the constant $c$ in the definition
of Assouad-Nagata dimension can be chosen the same for all $A_s$, $s\in S$.
Geometrically, the most common occurrence is when all $A_s$ are isometric.

Here is a generalization of \ref{DimZeroForMetricallyParallelFibers}:

\begin{Thm}\label{GeneralDimForMetricallyParallelFibers}
Suppose $f\colon (X,\rho_X)\to (Y,\rho_Y)$ satisfies $\rho_X(x,f^{-1}(y))\leq \rho_Y(f(x),y)$ for all $x\in X$ and $y\in Y$. If $f$ is Lipschitz and $$\dim_{AN}(f^{-1}(y))\leq k$$
uniformly with respect to $y\in f(X)$,
then $\dim_{AN}(X)\leq k+\dim_{AN}(Y)$.
\end{Thm}
\dokaz The proof consists of applying Theorem 7.1 in \cite{BDLM}.
That result says $\dim_{AN}(X)\leq \dim_{AN}(f)+\dim_{AN}(Y)$
if $f$ is Lipschitz.
We need to show $\dim_{AN}(f)\leq k$ which
means the existence of constants $a$ and $b$ such that
for any $r_X, R_Y>0$ and any subset $A$ of $X$ so that $f(A)$ is $R_Y$-bounded,
$A$ can be expressed as $A=A_0\cup\ldots \cup A_k$ with $r_X$-components
of each $A_i$ being $(a\cdot r_X+b\cdot R_Y)$-bounded.
\par Assume $f$ is $\lambda$-Lipschitz.
Pick $y\in f(A)$ and put $F=f^{-1}(y)$.
Cover $F$ with sets $F_0,\ldots,F_k$ such that
$(2R_Y+r_X)$-components of $F_i$ are $c\cdot(2R_Y+r_X)$-bounded.
Look at $A_i=B(F_i,R_Y)$ and notice its $r_X$-components
are $(c(2R_Y+r_X)+2R_Y)$-bounded.
As $A_i$ cover $A$ (use $\rho_X(x,F)\leq \rho_Y(f(x),f(F))< R_Y$ for all $x\in A$)
$a=c$ and $b=2c+2$ work.
\edokaz

\section{Dimension of asymptotic cones}

\begin{Prop}\label{ConesAndPnProperty}
Let $(X,\rho_X)$ be a metric space and $n\ge 0$.
If $(X,\rho_X)$ has $n$-dimensional Nagata property, then every asymptotic cone
over $X$ has $n$-dimensional Nagata property.
\end{Prop}
\dokaz Let's rephrase the $n$-dimensional Nagata property as follows:
given $r > 0$ and points $y_i,z_i$, $i\leq n+2$, such that $\rho_X(y_i,z_i) < r$
and $\rho_X(x,z_i) < r/2$ for some $x\in X$ and all $i\leq n+2$,
there is $i\ne j$ such that $\rho_X(y_i,y_j)< r$.
\par Suppose $(Y,\rho_Y)=\AC_\omega(X,c,d)$ does not have $n$-dimensional Nagata property
and there exist points $x\in Y$, $y_i,z_i\in Y$ for $i\leq n+2$
such that for some $r > \epsilon > 0$ the following inequalities hold:
\begin{itemize}
\item[a.] $\rho_Y(y_i,z_i) < r-\epsilon$ and $\rho_Y(x,z_i) < r/2-\epsilon/2$ for all $i\leq n+2$,
\item[b.] $\rho_Y(y_i,y_j)\ge r$ for all $i\ne j$.
\end{itemize}
For each $i$ pick $P_i\in\omega$ such that $\rho_X(y_i(n),z_i(n))< (r-\epsilon)\cdot d_n$
and $\rho_X(z_i(n),x(n))< (r/2-\epsilon/2)\cdot d_n$ for all $n\in P_i$.
Similarly, for all $i\ne j$, pick $P_{ij}\in \omega$ such that
$\rho_X(y_i(n),y_j(n))\ge (r-\epsilon)\cdot d_n$ for all $n\in P_{ij}$.
Let $P$ be the intersection of all $P_i$ and all $P_{ij}$.
Since $P\ne\emptyset$, we arrive at a contradiction: $(X,\rho_X)$
cannot have $n$-dimensional Nagata property.
\edokaz

\begin{Cor}\label{ConesAndUltrametrics}
If $(X,\rho_X)$ is ultrametric, then every asymptotic cone
$\AC_\omega(X,c,d)$ is ultrametric.
\end{Cor}
\dokaz
It follows from the well-known fact that being ultrametric
is the same as having $0$-dimensional Nagata property (easy exercise left to the reader).
\edokaz
\begin{Cor}\label{ConesAndANDim}
$$\dim(\AC_\omega(X,c,d))\leq \dim_{AN}(\AC_\omega(X,c,d))\leq \asdim_{AN}(X,\rho_X)$$ for any metric space $(X,\rho_X)$.
\end{Cor}
\dokaz Suppose $\asdim_{AN}(X,\rho_X)=n$. Notice $\AC_\omega(X,\rho_X,d)$
is isometric to $\AC_\omega(X,\rho_2,d)$, where $\rho_2(x,y)=\rho_X(x,y)$
if $\rho_X(x,y)\ge 1$ and $\rho_2(x,y)=1$ if $0 < \rho_X(x,y)< 1$.
Also, $\dim_{AN}(X,\rho_2)=n$.
By \ref{ANDimAndPnProperty} there is a metric space
$(Y,\rho_Y)$ bi-Lipschitz equivalent to $(X,\rho_2)$ such that
$(Y,\rho_Y^p)$ has $n$-dimensional Nagata property for some $0 < p < 1$.
Therefore $\AC_\omega(Y,\rho_Y,d)$ is bi-Lipschitz equivalent
to $\AC_\omega(X,\rho_X,d)$. Since the snowflaked version
$\AC_\omega(Y,\rho_Y^p,d^p)$ of $\AC_\omega(Y,\rho_Y,d)$ has $n$-dimensional Nagata property by \ref{ConesAndPnProperty},
we conclude $\dim_{AN}(\AC_\omega(X,\rho_X,d))\leq n$.
As $\dim(Y)\leq \dim_{AN}(Y)$ for all metric spaces (see \cite{LangSch}, Theorem 2.2),
the proof is completed.
\edokaz

J.Burillo \cite{Bur} showed that asymptotic cones of Baumslag-Solitar groups
$BS_{p,q}$ (groups with generators $t,x$ and the sole relation
$t^{-1}x^pt=x^q$ where $p\ne q$) have dimension $1$ by displaying
a map $f$ from the Cayley graph of $BS_{p,q}$ to an $\RR$-tree $T$
with metrically parallel fibers such that the fibers
of the induced maps on asymptotic cones are ultrametric. In view of
\ref{ConesAndANDim} and \ref{DimZeroForMetricallyParallelFibers}
one can improve his result by stating that asymptotic cones of
$BS_{p,q}$ are of Assouad-Nagata dimension $1$.

\begin{Question}
Is there a finitely generated group whose
asymptotic cones have bounded Assouad-Nagata dimension
but the group itself has infinite asymptotic Assouad-Nagata dimension?
\end{Question}

\begin{Question}
Is there $(X,\rho)$ of positive asymptotic Assouad-Nagata dimension
whose all asymptotic cones are of Assouad-Nagata dimension at most $0$?
\end{Question}
\section{Assouad-Nagata dimension of the discrete
Heisenberg group}

\par We do not know of a finitely presented group $G$ of finite asymptotic dimension
whose asymptotic Assouad-Nagata dimension is larger but finite.
Piotr Nowak \cite{Nowak} identified finitely generated groups $G_n$
of asymptotic dimension $n\ge 2$ and $\asdim_{AN} (G_n)=\infty$.
 A.N. Dranishnikov communicated to us that there was a suggestion by 
John Roe that perhaps the Heisenberg group
may provide an example where the dimensions differ.
However, that is not the case\footnote{When the paper was written we learned from Dranishnikov that he and P.Nowak proved this fact independently.}.

\begin{Prop}\label{Heisenberg}
The asymptotic dimension and Assouad-Nagata dimension of the discrete
Heisenberg group $H_3(\Z)$ both equal $3$.
\end{Prop}
\dokaz 
Recall (see \cite{GroAsymInvar}, p.51) the continuous Heisenberg group $H_3(\R)$ is the only
non-Abelian nilpotent three-dimensional group.
It contains, as a cocompact lattice, the discrete
Heisenberg group $H_3(\Z)$
generated by $a$, $b$, and $c$ with the relations $c=[a,b]$, $[a,c]=1$, and $[b,c]=1$.
Since $G=H_3(\Z)$ is of Hirsch length $3$, its asymptotic dimension is $3$
(see Theorem 3.5 in \cite{DS}).
To observe $\asdim_{AN}(G)\leq 3$ consider
the exact sequence $1\to \Z\to G\to G/\Z\to 1$, where $\Z$ is generated by $c$.
As $G/\Z=\Z\oplus\Z$, it suffices to show $(\Z,d)$
is quasi-symmetrically equivalent to $\Z$ equipped with the regular metric, where $d$ is a metric on $\Z$
induced from a word metric on $G$. Indeed, in view of \ref{LSThm} one gets $\asdim_{AN}(\Z,d)=1$ and applying
\ref{GeneralDimForMetricallyParallelFibers} one derives $3=\asdim(G)\leq \asdim_{AN}(G)\leq 3$.
\par 
As suggested in \cite{GroAsymInvar} on p.52
the metric $d$ is quasi-isometric to the square root of the regular metric on $\Z$
(the equality $[a^n,b^n]=c^{n^2}$ is cited). Since that statement is not totally obvious, let us provide
a detailed exposition of the following:
\par {\bf Claim 1}. Let $K > 16$ be sufficiently
large so that
$4(n+1)+K\cdot \sqrt{n}\leq K\cdot n$ for all $n\ge 2$
(equivalently, $K\ge \frac{4(n+1)}{n-\sqrt{n}}$ for all $n\ge 2$).
For every integer $k$ there is a word $w$ in $a$ and $b$ (including inverses) such that
$w=c^k$ in $G$ and $l(w)\leq K\sqrt{|k|}$,
where $l(w)$ is the length of $w$.
\par
{\bf Proof of Claim 1:}  We will use equality $[a^u,b^v]=c^{u\cdot v}$ which we leave to the reader to prove.
It suffices to consider $k$ positive and the proof is by induction on $k$.
If $k=n^2$ for some natural $n$, then we put $w=[a^n,b^n]$ whose length is $4\cdot n$.
For $k=2$ one takes $w=[a,b]^2$. Similarly, $k=3$ is handled.
If there is natural $n\ge 2$ such that $n^2 < k \leq n^2+n$,
then we pick a word $w'$ for $k'=k-n^2$ and we put
$w=[a^n,b^n]\cdot w'$. Notice $l(w)\leq 4n+l(w')\leq
4n+K\sqrt{k-n^2}\leq Kn\leq K\sqrt{k}$.
Otherwise there is $n\ge 2$ such that
$(n+1)^2-n\leq k < (n+1)^2$ in which case we put $k'=k-(n+1)^2$
and $w=[a^{n+1},b^{n+1}]\cdot w'$, where $w'=c^{k'}$ and $l(w')\leq K\sqrt{n}$.
\edokaz

\par {\bf Claim 2}. If $w$ is a word in $a$ and $b$ (including inverses) such that
$w=c^k$ in $G$ for some $k$, then $2\sqrt{|k|}\leq l(w)$.
\par
{\bf Proof of Claim 2:}  By contradiction.
Find $w$ of minimal length such that $w=c^k$ and $l(w) < 2 \sqrt{|k|}$.
Among such words (minimizing $l(w)$) choose a word maximizing
$|k|$ (that maximum cannot be infinity as there are only finitely many words
of a given length).
Notice $l(w) > 4$. Since $c$ is in the center of $G$, equality $x\cdot y=c^k$
implies $y\cdot x=c^k$. Therefore (via a cyclic permutation)
we may assume $w=a^{m(1)}\cdot b^{n(1)}\cdot \ldots a^{m(s)}\cdot b^{n(s)}$,
where all exponents are non-zero
and $l(w)=\sum (|m(i)|+|n(i)|)$. Notice $\sum m(i)=\sum n(i)=0$.
Also, we may assume $k > 0$.
\par The first observation is $m(i)\cdot n(i) > 0$ for all $i$.
Indeed, if $m(i)\cdot n(i) < 0$ for some $i$ we may assume $i=s$ without loss of generality
and observe $a^{m(s)}b^{n(s)}=b^{n(s)}a^{m(s)}c^{m(s)n(s)}$, so one can
find either a shorter word $w'$ or larger $k'=k-m(s)n(s)$ satisfying $w'=c^{k'}$,
a contradiction.
\par The same way the sign of $m(i)$ cannot be the same as of $m(i+1)$ for any $i$.
Indeed $b^{n(i)}a^{m(i+1)}=c^{-n(i)m(i+1)}a^{m(i+1)}b^{n(i)}$ and as above
we get a contradiction.
\par If one looks for the smallest absolute value of exponents $m(i)$ and $n(i)$,
then (up to cyclic permutation) one can express $w$ as
either $x(a^u\cdot b^u\cdot a^{-u})y$ 
or  $x(b^u\cdot a^{-u}\cdot b^{-u})y$ with $l(w)=l(x)+l(y)+3|u|$.
In the first case one has $a^u\cdot b^u\cdot a^{-u}=c^{u^2}b^u$,
so $x\cdot b^u\cdot y=c^{k-u^2}$.
By induction $l(w)-2|u|\ge 2\sqrt{|k-u^2|}$,
so $l(w)\ge 2|u|+2\sqrt{|k-u^2|}\ge 2\sqrt{k}$, a contradiction.
The second case is similar.
\edokaz
\par
Now $d(c^m,c^n)=d(c^{m-n},1)\ge 2\sqrt{|m-n|}$ by Claim 2
and $d(c^m,c^n)\leq K\sqrt{|m-n|}$ by Claim 1.
That means $(\Z,d)$ is bi-Lipschitz equivalent to $(\Z,\sqrt{|m-n|})$.
Since $(\Z,\sqrt{|m-n|})$ is quasi-symmetrically equivalent to $(\Z,|m-n|)$,
we are done.
\edokaz

It is not clear if using methods of \cite{Osin}
one can generalize \ref{Heisenberg} and prove
$\asdim(G)=\asdim_{AN}(G)$ for any virtually
nilpotent group $G$.

\begin{Question}
Is there a finitely generated group of Assouad-Nagata dimension
$1$ that is not finitely presented?
\end{Question}


\begin{thebibliography}{99}

\bibitem{Assouad}
P. Assouad, {\em Sur la distance de Nagata}, C. R. Acad. Sci. Paris Ser. I
Math. {\bf 294} (1982), no. 1, 31--34.

\bibitem{Behr}
J.~Behrstock, {\em Asymptotic geometry of the mapping class group and Teichmuller space},
preprint ArXiv:math.GT/0502367.

\bibitem{BehrMinsky}
J.~Behrstock and Y.N.~Minsky, {\em Dimension and rank for mapping class groups},
preprint ArXiv:math.GT/0512352.

\bibitem{BellDranish}
G. Bell and A. Dranishnikov, {\em 
 A Hurewicz-type theorem for asymptotic dimension and applications
to geometric group theory}, Preprint, math.GR/0407431 (2004).

\bibitem{BellDranish2}
G. Bell and A. Dranishnikov, {\em 
 Asymptotic dimension in Bedlewo}, Topology Proceedings (to appear).

\bibitem{BDHM}
N. Brodskiy, J. Dydak, J. Higes, A. Mitra, \emph{Nagata-Assouad dimension via Lipschitz extensions}, preprint ArXiv:math.MG/0601226.

\bibitem{BDLM}
N. Brodskiy, J. Dydak, M. Levin, A. Mitra, \emph{Hurewicz Theorem for Assouad-Nagata dimension}, preprint ArXiv:math.MG/0605416.

\bibitem{Bur}
J.Burillo, {\em Dimension and fundamental group of asymptotic cones},
Journal of the London Mathematical Society (1999), 59, 557--572.

\bibitem{DranPrList}
A.Dranishnikov, {\em private communication}.

   \bibitem{DS}
A.Dranishnikov and J.Smith,
{\em Asymptotic dimension of discrete groups}, Fund. Math. 189 (2006), 27--34.

\bibitem{DRS}
A.N.Dranishnikov and J.Smith, {\em On asymptotic Assouad-Nagata dimension}, 
preprint ArXiv: math.MG/0607143

\bibitem{Dran-Zar}
A.~Dranishnikov, M.~Zarichnyi, \emph{Universal spaces for
asymptotic dimension}, Topology Appl. {\bf 140} (2004), no.2-3,
203--225.

\bibitem{Drutu}
C.~Drutu, \emph{Quasi-isometry invariants and asymptotic cones}, Int. J. Algebra Comput. 12 (1 and 2) (2002),99--135.

\bibitem{DruSap}
C.~Drutu, M.~Sapir, \emph{Tree-graded spaces and asymptotic cones of groups}, Topology {\bf 40} (2005), 959--1058.


\bibitem{GroNilp}
M.~Gromov, {\em Groups of polynomial growth and expanding maps}, Publ. Math. IHES 53 (1981) 53--73.

\bibitem{GroAsymInvar}
M.~Gromov, {\em Asymptotic invariants for infinite groups}, in
Geometric Group Theory, vol. 2, 1--295, G. Niblo and M. Roller,
eds., Cambridge University Press, 1993.



\bibitem{Hei}
J.Heinonen, {\em Lectures on Analysis on Metric Spaces}, Springer Verlag, Universitext 2001.

\bibitem{LangSch}
U. Lang, T. Schlichenmaier, {\em Nagata dimension, quasisymmetric
embeddings, and Lipschitz extensions}, IMRN International
Mathematics Research Notices (2005), no. 58, 3625--3655.

\bibitem{Kap}
M.Kapovich, {\em Lectures on Geometric Group Theory},
preprint (as of September 28, 2005).

\bibitem{Nag}
J.Nagata, {\em Note on dimension theory of metric spaces}, Fund.Math.45(1958), 143--181.

 \bibitem{Nowak}
 P.W.Nowak, {\em On exactness and isoperimetric profiles of discrete groups}, preprint.

\bibitem{Osin}
D. V. Osin, {\em Subgroup distortions in
nilpotent groups}, Comm.Algebra 29 (2001), pp.5439--5463.

\bibitem{Pansu}
P.~Pansu, {\em Croissance des boules et des geodesiques fermees dans les nilvarietes}, Ergod. Th. Dynam. Syst. 3 (1983),
415--445.

\bibitem{Roe lectures}
J. Roe, {\em Lectures on coarse geometry}, University Lecture
Series, 31. American Mathematical Society, Providence, RI, 2003.

\bibitem{Dries}
L. van den Dries, A.J. Wilkie, {\em On Gromov's theorem concerning groups of polynomial growth and elementary logic}, J.
Algebra 89 (1984), 349--374.


\end{thebibliography}
\end{document}